\documentclass[11pt]{amsart}
\usepackage{amssymb}

\newtheorem{lemma}{Lemma}[section]

\newtheorem{theorem}[lemma]{Theorem}

\newtheorem{corollary}[lemma]{Corollary}
\newtheorem{definition}[lemma]{Definition}

\begin{document}

\vspace*{1cm}

\title{Stratified critical points on the real Milnor fibre and integral-geometric formulas}
\author{Nicolas Dutertre }
\address{Aix-Marseille Universit\'e, LATP,
39 rue F. Joliot-Curie,
13453 Marseille Cedex 13, France.}
\email{nicolas.dutertre@univ-amu.fr}

%N. Dutertre is supported by {\em Agence Nationale de la Recherche}
%(reference ANR-08-JCJC-0118-01)\\

\maketitle
\begin{center}
{\em Dedicated to professor David Trotman on his 60th birthday}
\end{center}

\begin{abstract} 
Let $(X,0) \subset (\mathbb{R}^n,0)$ be the germ of a closed subanalytic set and let $f$ and $g : (X,0) \rightarrow (\mathbb{R},0)$ be two subanalytic functions. Under some conditions, we relate the critical points  of $g$ on the real Milnor fibre $X \cap f^{-1}(\delta) \cap B_\epsilon$, $0 < \vert \delta \vert \ll \epsilon \ll 1$, to the topology of this fibre and other related subanalytic sets. As an application, when $g$ is a generic linear function, we obtain an ``asymptotic" Gauss-Bonnet formula for the real Milnor fibre of $f$. From this Gauss-Bonnet formula, we deduce  ``infinitesimal" linear kinematic formulas.

 \end{abstract}

\markboth{N. Dutertre}{Stratified critical points on the real Milnor fibre and integral-geometric formulas}

\section{Introduction}
Let $F=(f_1,\ldots,f_k) :  (\mathbb{C}^n,0) \rightarrow (\mathbb{C}^k,0)$, $2 \le k \le n$, be a complete intersection with isolated singularity. The L\^e-Greuel formula \cite{Greuel,Le} states that
$$\mu(F')+ \mu(F)= \hbox{dim}_{\mathbb{C}} \frac{\mathcal{O}_{\mathbb{C}^n,0}}{I}, $$
where $F' :(\mathbb{C}^n,0) \rightarrow (\mathbb{C}^{k-1},0)$ is the map with components $f_1,\ldots,f_{k-1}$,
$I$ is the ideal generated by $f_1,\ldots,f_{k-1}$ and the $(k \times k)$-minors
$\frac{\partial(f_1,\ldots,f_k)}{\partial(x_{i_1},\ldots,x_{i_k})}$ and $\mu(F)$ (resp. $\mu(F')$) is the Milnor number of $F$ (resp $F'$).
Hence the L\^e-Greuel formula gives an algebraic characterization of a topological data, namely the sum of two Milnor numbers.
However, since the right-hand side of the above equality is equal to the number of critical points of $f_k$,  counted with multiplicity, on the Milnor fibre of $F'$, the L\^e-Greuel formula can be also viewed as a topological characterization of this number of critical points. 

Many works have been devoted to the search of a real version of the L\^e-Greuel formula. Let us recall them briefly.
We consider an analytic map-germ $F =(f_1,\ldots,f_k) : (\mathbb{R}^n,0) \rightarrow (\mathbb{R}^k,0)$, $2 \le k \le n$, and we denote by $F'$ the map-germ $(f_1,\ldots,f_{k-1}): (\mathbb{R}^n,0) \rightarrow (\mathbb{R}^{k-1},0)$. Some authors investigated the following difference:
$$D_{\delta,\delta'}= \chi \big( F'^{-1}(\delta) \cap \{f_{k} \ge \delta'\} \cap B_\epsilon \big) - \chi \big( F'^{-1}(\delta) \cap \{f_{k} \le \delta'\} \cap B_\epsilon \big),$$
where $(\delta,\delta')$ is a regular value of $F$ such that $0 \le \vert \delta' \vert \ll \vert \delta \vert \ll \epsilon$.

In \cite{DutertreHokkaido}, we proved that
$$D_{\delta,\delta'} \equiv \hbox{dim}_{\mathbb{R}} \frac{\mathcal{O}_{\mathbb{R}^n,0}}{I} \bmod 2,$$
where $\mathcal{O}_{\mathbb{R}^n,0}$ is the ring of analytic function-germs at the origin and $I$ is the ideal generated by $f_1,\ldots,f_{k-1}$ and all the $k\times k$ minors $\frac{\partial (f_k,f_1,\ldots,f_{k-1})}{\partial(x_{i_1},\ldots,x_{i_k})}$. This is only a mod 2 relation and we may ask if it is possible to get a more precise relation.

When $k=n$ and $f_k=x_1^2+\cdots+x_n^2$, according to Aoki et al. (\cite{AFN1}, \cite{AFS}), $D_{\delta,0}= \chi \big(F'^{-1}(\delta) \cap B_\varepsilon \big)= 2 \hbox{deg}_0 H$ and $2\hbox{deg}_0 H$ is the number of semi-branches of $F'^{-1}(0),$ where
$$H=(\frac{\partial (f_n,f_1,\ldots,f_{n-1})}{\partial(x_{1},\ldots,x_{n})},f_1,\ldots,f_{n-1}).$$
They proved a similar formula in the case $f_k=x_n$ in \cite{AFN2} and Szafraniec generalized all these results to any $f_k$ in \cite{Szafraniec}.

When $k=2$ and $f_2=x_1$, Fukui \cite{FukuiHawaii} stated that $$D_{\delta,0}=-\hbox{sign}(-\delta)^n \hbox{deg}_{0} H,$$ where $H=(f_1,\frac{\partial f_1}{\partial x_2},\ldots,\frac{\partial f_1}{\partial x_n})$. Several generalizations of Fukui's formula are given in \cite{FukuiTopology}, \cite{DutertreKodai1}, \cite{FukuiKhovanskii} and \cite{DutertreKodai2}. 

In all these papers, the general idea is to count algebraically the critical points of a Morse perturbation of $f_k$ on $F'^{-1}(\delta) \cap B_\epsilon$ and to express this sum in two ways: as a difference of Euler characteristics and as a topological degree. Using the Eisenbud-Levine formula \cite{EisenbudLevine}, this latter degree can be expressed as a signature of a quadratic form and so, we obtain an algebraic expression for $D_{\delta,\delta'}$. 

In this paper, we give a real and stratified version of the L\^e-Greuel formula. We restrict ourselves to the topological aspect and relate a sum of indices of critical points on a real Milnor fibre to some Euler characteristics (this is also the point of view adopted in \cite{CisnerosGrulhaSeade}).
More precisely, we consider  a germ of a closed subanalytic set $(X,0) \subset (\mathbb{R}^n,0)$  and a subanalytic function $f : (X,0) \rightarrow (\mathbb{R},0)$. We assume that $X$ is contained in a open set $U$ of $\mathbb{R}^n$ and that $f$ is the restriction to $X$ of a $C^2$-subanalytic function $F : U \rightarrow \mathbb{R}$.  We denote by $X^f$ the set $X \cap f^{-1}(0)$ and we equip $X$ with a Thom stratification adapted to $X^f$.  If $0 <\vert \delta \vert \ll \epsilon \ll 1$ then the real Milnor fibre of $f$ is defined by 
$$ M_f^{\delta,\epsilon}=f^{-1}(\delta) \cap X \cap B_\epsilon.$$
We consider another  subanalytic function $g : (X,0) \rightarrow (\mathbb{R},0)$ and we assume that it is the restriction to $X$ of a $C^2$-subanalytic function $G: U \rightarrow \mathbb{R}$. We denote by $X^g$ the intersection $X \cap g^{-1}(0)$. Under two conditions on $g$, we  study the topological behaviour of $g_{\vert M_f^{\delta,\epsilon}}$. 

We recall that 
if $Z \subset \mathbb{R}^n$ is a closed subanalytic set, equipped with a Whitney stratification and $p \in Z$ is an isolated critical point of a subanalytic function $\phi : Z \rightarrow \mathbb{R}$, restriction to $Z$ of a $C^2$-subanalytic function $\Phi$, then the index of $\phi$ at $p$ is defined as follows:
$${\rm ind}(\phi,Z,p)= 1 - \chi \big( Z \cap \{ \phi = \phi(p)-\eta \} \cap B_\epsilon(p) \big),$$
where $0< \eta \ll \epsilon \ll 1$ and $B_\epsilon(p)$ is the closed ball of radius $\epsilon$ centered at $p$. 
Let $p_1^{\delta,\epsilon},\ldots,p_r^{\delta,\epsilon} $ be the critical points of $g$ on $X \cap f^{-1}(\delta) \cap \mathring{B_\epsilon}$, where $\mathring{B_\epsilon}$ denotes the open ball of radius $\epsilon$. We set 
$$I(\delta, \epsilon, g) = \sum_{i=1}^r {\rm ind}(g, X \cap f^{-1}(\delta), p_i ^{\delta,\epsilon}),$$
$$I(\delta, \epsilon, -g) = \sum_{i=1}^r {\rm ind}(-g,X \cap f^{-1}(\delta), p_i ^{\delta,\epsilon}).$$
Our main theorem (Theorem \ref{MainTh}) is the following:
$$I(\delta,\epsilon,g)+I(\delta,\epsilon,-g) = 2 \chi(M_f^{\delta,\epsilon}) 
 - \chi (X \cap f^{-1}(\delta) \cap S_\epsilon)- \chi ( X^g \cap f^{-1}(\delta) \cap S_\epsilon). $$
As a corollary (Corollary \ref{MainCor}), when $f : (X,0) \rightarrow (\mathbb{R},0)$ has an isolated stratified critical point at $0$, we obtain that
$$I(\delta,\epsilon,g) + I(\delta,\epsilon,-g) = 2 \chi ( M_f^{\delta,\epsilon}) - \chi ({\rm Lk}(X^f)) - \chi ({\rm Lk} (X^f \cap X^g)),$$
where Lk$(-)$ denotes the link at the origin.

Then we apply these results when $g$ is a generic linear form to get an asymptotic Gauss-Bonnet formula for $M_f^{\delta,\epsilon}$ (Theorem \ref{GB}). 
In the last section, we use this asymptotic Gauss-Bonnet formula to prove infinitesimal linear kinematic formulas for closed subanalytic germs (Theorem \ref{KinForm}), that generalize the Cauchy-Crofton formula for the density due to Comte \cite{Comte}. 

The paper is organized as follows. In Section 2, we prove several lemmas about critical points on the link of a subanalytic set. Section 3 contains real stratified versions of the L\^e-Greuel formula.
In Section 4, we establish the asymptotic Gauss-Bonnet formula and in Section 5, the infinitesimal linear kinematic formulas.

%We will use the following notations: for $p \in \mathbb{R}^n$ and $\varepsilon >0$, $B_\varepsilon^n(p)$ is the ball of radius $\varepsilon$ centered at $p$ and $S_\varepsilon^{n-1}(p)$ the sphere of radius $\varepsilon$ centered at $p$. If $p=0$, we simply set $B_\varepsilon^n $ and $S_\varepsilon^{n-1}$ and if $p=0$ and $\varepsilon=1$ we use the standard notations $B^n$ and $S^{n-1}$. If $\mathcal{E}$ is a subset of $\mathbb{R}^n$ then $\mathring{\mathcal{E}}$ denotes its topological interior.
%Before stating our main results, we need some notations:
%\begin{itemize}
%\item for $k \in \{0,\ldots,n \}$, $G_{n}^k$ is the Grassmann manifold of $k$-dimensional linear
%subspaces in $\mathbb{R}^{n}$ and $g_n^k$ is its volume,
%\item for $k \in \mathbb{N}$, $b_k$ is the volume of the $k$-dimensional unit ball and $s_k$ is the
%volume of the $k$-dimensional unit sphere,
%\item  for $R>0$, $B_R$ (resp. $S_R$) will denote the ball (resp. the sphere) centered at the origin of radius $R$ in $\mathbb{R}^n$.
%If $R=1$, we will write $B^n$ (resp. $S^{n-1}$).
%\end{itemize}

The author is  grateful to Vincent Grandjean for a very useful discussion on generic distance functions. 

The author is partially supported by the program 
\begin{center}``Cat\'edras L\'evi-Strauss$-$USP/French Embassy, no. 2012.1.62.55.7". \end{center}
This paper was written while the author was visiting the {\em Instituto de Ci\^{e}ncias Matem\'{a}ticas e de Computa\c{c}\~{a}o, Universidade de S\~{a}o Paulo - Campus de S\~{a}o Carlos}. He thanks this institution, especially Raimundo Ara\'ujo dos Santos and Nivaldo Grulha, for the hospitality.

\section{Lemmas on critical points on the link of a stratum}\label{qualquer}
In this section, we study the behaviour of the critical points of a $C^2$-subanalytic function on the link of stratum that contains $0$ in its closure, for  a generic choice of the distance function to the origin.

Let $Y \subset \mathbb{R}^n$ be a $C^2$-subanalytic set such that $0$ belongs to its closure $\overline{Y}$.  Let $\theta : \mathbb{R}^n \rightarrow \mathbb{R}$ be a $C^2$-subanalytic function such that $\theta (0)=0$. We will first study the behaviour of the critical points of $\theta_{\vert Y} : Y \rightarrow \mathbb{R}$ in the neighborhood  of $0$, and then the behaviour of the critical points of the restriction of $\theta$ to the link of $0$ in $Y$.
\begin{lemma}\label{CrPts1}
The critical points of $\theta_{\vert Y}$ lie in $\{\theta = 0 \}$ in a neighborhood of $0$.
\end{lemma}
\proof By the Curve Selection Lemma, we can assume that there is a $C^1$-subanalytic curve $\gamma : [0, \nu [ \rightarrow \overline{Y}$ such that $\gamma(0)=0$ and $\gamma (t)$ is a critical point of $\theta_{\vert Y}$ for $t  \in ]0,\nu [$. Therefore, we have 
$$(\theta \circ \gamma)'(t)= \langle \nabla \theta_{\vert Y} (\gamma(t)), \gamma'(t) \rangle =0,$$
since $\gamma'(t)$ is tangent to $Y$ at $\gamma(t)$. This implies that $\theta \circ \gamma (t)= \theta \circ \gamma (0)=0$. \endproof

Let $\rho : \mathbb{R}^n \rightarrow \mathbb{R}$ be another $C^2$-subanalytic function such that $\rho^{-1}(a)$ intersects $Y$ transversally. Then the set $Y \cap \{ \rho \le a \}$  is a manifold with boundary. Let $p$ be a critical point of $\theta_{\vert Y \cap \{\rho \le a \}}$ which lies in $Y \cap \{ \rho = a \}$ and which is not a critical point of $\theta_{\vert Y}$. This implies that 
$$\nabla \theta_{\vert Y}(p) = \lambda (p) \nabla \rho_{\vert Y} (p),$$
with $\lambda (p) \not= 0$. 

\begin{definition}
We say that $p \in Y \cap \{ \rho = a \}$ is an outwards-pointing (resp. inwards-pointing) critical point of $\theta_{\vert Y \cap \{ \rho \le a \}}$ if $\lambda (p) >0$ (resp. $\lambda (p) <0$).
\end{definition}
Now let us assume that $\rho : \mathbb{R}^n \rightarrow \mathbb{R}$ is a distance function to the origin which means that $\rho \ge 0$ and $\rho^{-1}(0) = \{ 0 \}$ in a neighborhood of $0$. By Lemma \ref{CrPts1}, we know that for $\epsilon >0$ small enough, the level $\rho^{-1}(\epsilon)$ intersects $Y$ transversally. Let $p^\epsilon$ be a critical point of $\theta_{\vert Y \cap \rho^{-1}(\epsilon)}$ such that $\theta (p^\epsilon) \not= 0$. 
This means that there exists $\lambda (p^\epsilon)$ such that 
$$\nabla \theta_{\vert Y} (p^\epsilon) = \lambda (p^\epsilon) \nabla \rho_{\vert Y}(p^\epsilon).$$
Note that $\lambda (p^\epsilon) \not= 0$ because $\nabla \theta_{\vert Y} (p^\epsilon) \not= 0$ for $\theta(p^\epsilon)  \not= 0$. 
\begin{lemma}\label{CrPts2}
The point $p^\epsilon$ is an outwards-pointing (resp. inwards-pointing) for $\theta_{\vert Y \cap \{ \rho \le \epsilon \}}$ if and only if $\theta(p^\epsilon) >0$ (resp. $\theta (p^\epsilon) <0$).
\end{lemma}
\proof Let us assume that $\lambda (p^\epsilon) >0$. By the Curve Selection Lemma, there exists a $C^1$-subanalytic curve $\gamma : [0,\nu [ \rightarrow \overline{Y}$ passing through $p^\epsilon$ such that $\gamma (0)=0$ and for $t \not= 0$, $\gamma (t)$ is a critical point of $\theta_{\vert Y \cap \{ \rho = \rho(\gamma(t)) \}}$ with $\lambda (\gamma(t)) >0$. Therefore we have 
$$(\theta \circ \gamma)'(t)= \langle \nabla \theta_{\vert Y} (\gamma(t)), \gamma'(t) \rangle = \lambda(\gamma (t)) \langle \nabla \rho_{\vert Y} (\gamma (t)), \gamma'(t) \rangle.$$
But $(\rho \circ \gamma)' >0$ for otherwise $(\rho \circ \gamma)' \le 0$ and $\rho \circ \gamma$ would be decreasing. Since $\rho (\gamma (t))$ tends to $0$ as $t$ tends to $0$, this would imply that $\rho \circ \gamma (t) \le 0$, which is impossible. We can conclude that $(\theta \circ \gamma)' >0$ and that $\theta \circ \gamma$ is strictly increasing. Since $\theta \circ \gamma (t)$ tends to $0$ as $t$ tends to $0$, we see that $\theta \circ \gamma (t) >0$ for $t >0$. Similarly if $\lambda (p^\epsilon) <0$ then $\theta (p^\epsilon) <0$. \endproof

Now we will study these critical points for a generic choice of the distance function. We denote by Sym$(\mathbb{R}^n)$ the set of symmetric $n \times n$-matrices with real entries, by Sym$^*(\mathbb{R}^n)$ the open dense subset of such matrices with non-zero determinant and by Sym$^{*,+}(\mathbb{R}^n)$ the open subset of these invertible matrices that are positive definite or negative definite. Note that these sets are semi-algebraic. 
For each $A \in {\rm Sym}^{*,+} (\mathbb{R}^n)$, we denote by $\rho_A$ the following quadratic form:
$$\rho_A(x) = \langle A x, x \rangle.$$

We denote by $\Gamma_{\theta,A}^Y$ the following subanalytic polar set:
$$\Gamma_{\theta,A}^Y = \left\{ x \in Y \ \vert \ {\rm rank}\left[ \nabla \theta_{\vert Y} (x), \nabla {\rho_A}_{\vert Y}(x) \right] < 2 \right\},$$
and by $\Sigma^Y_\theta$ the set of critical points of $\theta_{\vert Y}$. Note that $\Sigma^Y_\theta \subset \{ \theta=0 \}$ by Lemma \ref{CrPts1}.
\begin{lemma}\label{CrPts3}
For almost all $A$ in ${\rm Sym}^{+,*} (\mathbb{R}^n)$, $\Gamma_{\theta,A}^Y \setminus ( \Sigma^Y_\theta \cup \{0\} )$ is a $C^1$-subanalytic curve (possible empty) in a neighborhood of $0$.
\end{lemma}
\proof We can assume that dim $Y >1$. Let 
$$\displaylines{
\quad Z = \Big\{ (x,A) \in \mathbb{R}^n \times {\rm Sym}^{+,*}(\mathbb{R}^n) \ \vert \ x \in Y \setminus ( \Sigma^Y_\theta \cup \{0\} ) \hfill \cr
\hfill \hbox{ and rank}\left[ \nabla \theta_{\vert Y} (x), \nabla {\rho_A}_{\vert Y}(x) \right] < 2 \Big\}. \quad \cr
}$$
Let $(y,B)$ be a point in $Z$. We can suppose that around $y$, $Y$ is defined by the vanishing of $k$ subanalytic functions $f_1,\ldots,f_k$ of class $C^2$. Hence in a neighborhood of $(y,B)$, $Z$ is defined be the vanishing of $f_1,\ldots,f_k$ and the minors 
$$\frac{\partial (f_1,\ldots,f_k,\theta,\rho_A)}{\partial (x_{i_1},\ldots,x_{i_{k+2}})}.$$
Furthermore, since $y$ does not belong to $\Sigma^Y_\theta$, we can assume that 
$$\frac{\partial (f_1,\ldots,f_k,\theta)}{\partial(x_1,\ldots,x_k,x_{k+1})} \not= 0,$$
in a neighborhood of $y$. Therefore $Z$ is locally defined by $f_1=\cdots=f_k=0$ and 
$$\frac{\partial (f_1,\ldots,f_k,\theta,\rho_A)}{\partial (x_1,\ldots,x_{k+1},x_{k+2})}= \cdots = \frac{\partial (f_1,\ldots,f_k,\theta,\rho_A)}{\partial (x_1,\ldots,x_{k+1},x_n)}=0.$$
Let us write $M=\frac{\partial(f_1,\ldots,f_k,\theta)}{\partial(x_1,\ldots,x_k,x_{k+1})}$ and for $i \in \{k+2,\ldots,n \}$, $m_i= \frac{\partial (f_1,\ldots,f_k,\theta,\rho_A)}{\partial (x_1,\ldots,x_{k+1},x_i)}$. If $A=[a_{ij}]$ then
$$\rho_A(x)= \sum_{i=1}^n a_{ii}x_i^2 + 2 \sum_{i\not= j} a_{ij}x_ix_j,$$
 and so $\frac{\partial \rho_A}{\partial x_i}(x)= 2 \sum_{j=1}^n a_{ij} x_j.$
For $i \in \{k+1,\ldots,n \}$ and $j \in \{1,\ldots,n \}$, we have
 $$\frac{\partial m_i}{\partial a_{ij}} = 2 x_j M.$$
 Since $y  \not= 0$, one of the $x_j$'s does not vanish in the neighborhood of $y$ and we can conclude that the rank of 
 $$\left[ \nabla f_1(x), \ldots, \nabla f_k (x), \nabla m_{k+2}(x,A),\ldots,\nabla m_n(x,A) \right]$$ 
 is $n-1$ and that $Z$ is a $C^1$-subanalytic manifold of dimension $\frac{n(n+1)}{2}+1$. 
Now let us consider the projection $\pi_2 : Z \rightarrow {\rm Sym}^{+,*}(\mathbb{R}^n)$, $(x,A) \mapsto A$. Bertini-Sard's theorem implies that the set $D_{\pi_2}$ of critical values of $\pi_2$ is a subanalytic set of dimension strictly less than $\frac{n(n+1)}{2}$. Hence, for all $A \notin D_{\pi_2}$, $\pi_2^{-1}(A)$ is a $C^1$-subanalytic curve (possibly empty). But this set is exactly $\Gamma_{\theta,A}^Y \setminus ( \Sigma^Y_\theta \cup \{0\} )$. \endproof

Let $R \subset Y$ be a subanalytic set of dimension strictly less than ${\rm dim}\ Y$. We will need the following lemma.
\begin{lemma}\label{CrPts4}
For almost all $A$ in ${\rm Sym}^{+,*} (\mathbb{R}^n)$, $\Gamma_{\theta,A}^Y \setminus ( \Sigma^Y_\theta \cup \{0\} ) \cap R$ is a subanalytic set of dimension at most $0$ in  a neighborhood of $0$.
\end{lemma}
\proof Let us put $l={\rm dim}\ Y$. Since $R$ admits a locally finite subanalytic stratification, we can assume that $R$ is a  $C^2$-subanalytic manifold of dimension $d$  with $d<l$. Let $W$ be the following subanalytic set:
$$\displaylines{
\quad W = \Big\{ (x,A) \in \mathbb{R}^n \times {\rm Sym}^{+,*}(\mathbb{R}^n) \ \vert \ x \in R \setminus ( \Sigma^Y_\theta \cup \{0\} ) \hfill \cr
\hfill \hbox{ and rank}\left[ \nabla \theta_{\vert Y} (x), \nabla {\rho_A}_{\vert Y}(x) \right] < 2 \Big\}. \quad \cr
}$$
Using the same method as in the previous lemma, we can prove that $W$ is a $C^1$-subanalytic manifold of dimension $\frac{n(n+1)}{2}+1+d-l$ and conclude, remarking that $d-l \le -1$. \endproof

Now we introduce a new $C^2$-subanalytic function $\beta : \mathbb{R}^n \rightarrow \mathbb{R}$ such that $\beta (0)=0$. We denote by  $\Gamma_{\theta,\beta,A}^Y$ the following subanalytic polar set:
$$\Gamma_{\theta,\beta,A}^Y = \left\{ x \in Y \ \vert \ {\rm rank}\left[ \nabla \theta_{\vert Y} (x), \nabla \beta_{\vert Y} (x), \nabla {\rho_A}_{\vert Y} (x) \right] < 3 \right\},$$
and by $\Gamma_{\theta,\beta}^Y $ the following subanalytic polar set:
$$\Gamma_{\theta,\beta}^Y = \left\{ x \in Y \ \vert \ {\rm rank}\left[ \nabla \theta_{\vert Y} (x), \nabla \beta_{\vert Y} (x) \right] < 2 \right\}.$$
\begin{lemma}\label{CrPts5}
For almost all $A$ in ${\rm Sym}^{+,*} (\mathbb{R}^n)$, $\Gamma_{\theta,\beta,A}^Y  \setminus (\Gamma_{\theta,\beta}^Y \cup \{0\} )$ is a $C^1$-subanalytic set of dimension at most $2$ (possibly empty) in  a neighborhood of $0$.
\end{lemma}
\proof We can assume that dim $Y >2$. 
Let 
$$\displaylines{
\quad Z = \Big\{ (x,A) \in \mathbb{R}^n \times {\rm Sym}^{+,*}(\mathbb{R}^n) \ \vert \ x \in Y, {\rm rank} \left[\nabla \theta_{\vert Y}(x),\nabla \beta_{\vert Y}(x) \right] = 2 \hfill \cr
\hfill \hbox{ and rank} \left[\nabla \theta_{\vert Y}(x),\nabla \beta_{\vert Y}(x), \nabla {\rho_A}_{\vert Y}(x) \right] <3 \Big\}. \quad \cr
}$$
Let $(y,B)$ be a point in $Z$. We can suppose that around $y$, $Y$ is defined by the vanishing of $k$ subanalytic functions $f_1,\ldots,f_k$ of class $C^2$. Hence in a neighborhood of $(y,B)$, $Z$ is defined by the vanishing of $f_1,\ldots,f_k$ and the minors 
$$\frac{\partial(f_1,\ldots,f_k,\theta,\beta,\rho_A)}{\partial (x_{i_1},\ldots,x_{i_{k+3}})}.$$
Since $y$ does not belong to $\Gamma_{\theta,\beta}^Y$, we can assume that
$$\frac{\partial(f_1,\ldots,f_k,\theta,\beta)}{\partial(x_1,\ldots,x_k,x_{k+1},x_{k+2})} \not= 0,$$
in a neighborhood of $y$. Therefore $Z$ is locally defined by $f_1,\ldots,f_k=0$ and 
$$\frac{\partial (f_1,\ldots,f_k, \theta, \beta, \rho_A)}{\partial (x_1,\ldots,x_{k+2},x_{k+3})}= \cdots= \frac{\partial (f_1,\ldots,f_k, \theta, \beta, \rho_A)}{\partial (x_1,\ldots,x_{k+2},x_{n})}=0.$$
It is clear that we can apply the same method as Lemma \ref{CrPts3} to get the result. \endproof

\section{L\^e-Greuel type formula}
In this section, we prove the L\^e-Greuel type formula announced in the introduction.

Let $(X,0) \subset (\mathbb{R}^n,0)$ be the germ of a closed subanalytic set and let $f : (X,0) \rightarrow (\mathbb{R},0)$ be a subanalytic function. We assume that $X$ is contained in a open set $U$ of $\mathbb{R}^n$ and that $f$ is the restriction to $X$ of a $C^2$-subanalytic function $F : U \rightarrow \mathbb{R}$.  We denote by $X^f$ the set $X \cap f^{-1}(0)$ and by \cite{Bekka}, we can equip $X$ with a Thom stratification $\mathcal{V}= \{ V_\alpha \}_{\alpha \in A}$ adapted to $X^f$. This means that $\{ V_\alpha \in \mathcal{V} \ \vert \ V_\alpha \nsubseteq X^f \}$ is a Whitney stratification of $X \setminus X^f$ and that for any pair of strata $(V_\alpha,V_\beta)$ with $V_\alpha \nsubseteq X^f$ and $V_\beta \subset X^f$, the Thom condition is satisfied.

Let us denote by $\Sigma_\mathcal{V} f$ the critical locus of $f$. It is the union of the critical loci of $f$ restricted to each stratum, i.e. $\Sigma_\mathcal{V} f= \cup_\alpha \Sigma (f_{\vert V_\alpha})$, where $\Sigma( f_{\vert V_\alpha})$ is the critical set of $f _{\vert V_\alpha} :V_\alpha \rightarrow \mathbb{R}$.  Since $\Sigma_{\mathcal{V}} f \subset f^{-1}(0)$ (see Lemma \ref{CrPts1}), the fibre $f^{-1}(\delta)$ intersects the strata $V_\alpha$'s, $V_\alpha \nsubseteq X^f$, transversally if $\delta $ is sufficiently small. Hence it is Whitney stratified with the induced stratification $\{f^{-1}( \delta) \cap V_\alpha \  \vert \ V_\alpha \nsubseteq X^f \}$.

By Lemma \ref{CrPts1}, we know that if $\epsilon >0$ is sufficiently small then the sphere $S_\epsilon$ intersects $X^f$ transversally. By the Thom condition, this implies that there exists $\delta(\epsilon) >0$ such that for each $\delta$ with $0 < \vert \delta \vert \le \delta(\epsilon)$, the sphere $S_\epsilon$ intersects the fibre $f^{-1}(\delta)$ transversally as well. Hence the set $f^{-1}(\delta) \cap B_\epsilon$ is a Whitney stratified set equipped with the following stratification:
$$\{ f^{-1}(\delta) \cap V_\alpha \cap \mathring{B_\epsilon}, f^{-1}(\delta) \cap V_\alpha \cap S_\epsilon \ \vert \ V_\alpha \nsubseteq X^f \}.$$

\begin{definition}
We call the set $f^{-1}(\delta ) \cap X \cap B_\epsilon$, where $0 < \vert \delta \vert \ll \epsilon \ll 1$, a real Milnor fibre of $f$.
\end{definition}
We will use the following notation: $M_f^{\delta,\epsilon}=f^{-1}(\delta) \cap X \cap B_\epsilon$.

Now we consider another subanalytic function $g : (X,0) \rightarrow (\mathbb{R},0)$ and we assume that it is the restriction to $X$ of a $C^2$-subanalytic function $G: U \rightarrow \mathbb{R}$. We denote by $X^g$ the intersection $X \cap g^{-1}(0)$. Under some restrictions on $g$, we will study the topological behaviour of $g_{\vert M_f^{\delta,\epsilon}}$.

First we assume that $g$ satisfies the following Condition (A):
\begin{itemize}
\item Condition (A): $g : (X,0) \rightarrow (\mathbb{R},0)$ has an isolated stratified critical point at $0$.
\end{itemize}
This means that for each strata $V_\alpha$ of $\mathcal{V}$, $g : V_\alpha \setminus \{0\} \rightarrow \mathbb{R}$ is a submersion in a neighborhood of the origin. 

In order to give the second assumption on $g$, we need to introduce some polar sets.
Let $V_\alpha$ be a stratum of $\mathcal{V}$ not contained in $X^f$. Let $\Gamma_{f,g}^{V_\alpha}$ be the following set:
$$ \Gamma_{f,g}^{V_\alpha} = \left\{ x \in V_\alpha \ \vert \ {\rm rank}[ \nabla f_{\vert V_\alpha}(x), \nabla g_{\vert V_\alpha}(x) ] < 2 \right\},$$
and let $\Gamma_{f,g}$ be the union $\cup \Gamma_{f,g}^{V_\alpha}$ where $V_\alpha \nsubseteq X^f$. We call $\Gamma_{f,g}$ the relative polar set of $f$ and $g$ with respect to the stratification $\mathcal{V}$. 
We will assume that $g$ satifies the following Condition (B):
\begin{itemize}
\item Condition (B): the relative polar set $\Gamma_{f,g}$ is a 1-dimensional $C^1$-subanalytic set (possibly empty) in a neighborhood of the origin.
\end{itemize}
Note that Condition (B) implies that $\overline{\Gamma_{f,g}} \cap X^f \subset \{0\}$ in a neighborhood of the origin because the frontiers of the $\Gamma_{f,g}^{V_\alpha}$'s are $0$-dimensional. 

From Condition (A) and Condition (B), we can deduce the following result. 
\begin{lemma}\label{Lemma1}
We have $\overline{\Gamma_{f,g} } \cap X^g \subset \{0\}$ in a neighborhood of the origin.
\end{lemma}
\proof If it is not the case then there is a $C^1$-subanalytic curve $\gamma : [0,\nu [ \rightarrow \Gamma_{f,g} \cap X^g$ such that $\gamma(0)=0$ and $\gamma(]0,\nu [) \subset X^g \setminus \{0\}$. We can also assume that $\gamma ( ]0, \nu [)$ is contained in a stratum $V$. For $t \in ]0,\nu [$, we have
$$0=(g \circ \gamma)'(t)= \langle \nabla g_{\vert V} (\gamma(t)), \gamma'(t) \rangle.$$
Since $\gamma(t)$ belongs to $\Gamma_{f,g}$ and $\nabla g_{\vert V}(\gamma(t))$ does not vanish for $g : (X,0) \rightarrow (\mathbb{R},0)$ has an isolated stratified critical point at $0$, we can conclude that $\langle \nabla f_{\vert V}(\gamma(t)), \gamma'(t) \rangle=0$ and that $(f \circ \gamma)'(t)=0$ for all $t \in ]0,\nu [$. Therefore $f\circ \gamma \equiv 0$ because $f(0)=0$ and $\gamma ([0,\nu [)$ is included in $X^f$. This is impossible by  the above remark. \endproof

Let $\mathcal{B}_1,\ldots,\mathcal{B}_l$ be the connected components of $\Gamma_{f,g} $, i.e. $\Gamma_{f,g} = \sqcup_{i=1}^l \mathcal{B}_i$. Each $\mathcal{B}_i$ is a $C^1$-subanalytic curve along which $f$ is strictly increasing or decreasing and the intersection points of the $\mathcal{B}_i$'s with the fibre $M_f^{\delta,\epsilon}$ are exactly the critical points (in the stratified sense) of $g$ on $X \cap f^{-1}(\delta) \cap \mathring{B_\epsilon}$. Let us write 
$$M_f^{\delta,\epsilon} \cap \sqcup_{i=1}^l \mathcal{B}_i = \{ p_1^{\delta,\epsilon},\ldots,p_r^{\delta,\epsilon} \} .$$
Note that  $r \le l$.  

Let us recall now the definition of the index of an isolated stratified critical point.
\begin{definition}
Let $Z \subset \mathbb{R}^n$ be a closed subanalytic set, equipped with a Whitney stratification. Let $p \in Z$ be an isolated critical point of a subanalytic function $\phi : Z \rightarrow \mathbb{R}$, which is the restriction to $Z$ of a $C^2$-subanalytic function $\Phi$. We define the index of $\phi$ at $p$ as follows :
$${\rm ind}(\phi,Z,p)= 1 - \chi \big( Z \cap \{ \phi = \phi(p)-\eta \} \cap B_\epsilon(p) \big),$$
where $0< \eta \ll \epsilon \ll 1$ and $B_\epsilon(p)$ is the closed ball of radius $\epsilon$ centered at $p$.
\end{definition}

Our aim is to give a topological interpretation to the following sum:
$$ \sum_{i=1}^r {\rm ind}(g,X \cap f^{-1}(\delta),p_i^{\delta,\epsilon}) + {\rm ind}(-g,X \cap f^{-1}(\delta),p_i^{\delta,\epsilon}).$$
For this, we will apply stratified Morse theory to $g_{\vert M_f^{\delta,\epsilon}} $. Note that the points $p_i$'s  are not the only critical points of $g_{\vert M_f^{\delta,\epsilon}}$ and other critical points can occur on the ``boundary" $M_f^{\delta,\epsilon} \cap S_\epsilon$.

The next step is to study the behaviour of these ``boundary" critical points for a generic choice of the distance function to the origin.
Let $\rho : \mathbb{R}^n \rightarrow \mathbb{R}$ be a $C^2$-subanalytic function which is a distance function to the origin. We denote by $\tilde{S}_\epsilon $ the level $\rho^{-1}(\epsilon)$ and by $\tilde{B}_\epsilon$ the set $\{ \rho \le \epsilon\}$. We will focus on the critical points of $g_{\vert X^f \cap \tilde{S}_\epsilon}$ and $g_{\vert X \cap f^{-1}( \delta) \cap \tilde{S}_\epsilon}$, with $0<\vert \delta \vert \ll \epsilon \ll 1$. 

For each stratum $V$ of $X^f$, let 
$$\Gamma_{g,\rho}^V = \left\{ x \in V \ \vert \ {\rm rank} [ \nabla g_{\vert V}(x), \nabla \rho_{\vert V}(x) ] < 2 \right\},$$
and let $\Gamma_{g,\rho}^{X^f} = \cup_{V \subset X^f} \Gamma_{g,\rho}^V$.  By Lemma \ref{CrPts3} and the fact that $g:(X^f,0) \rightarrow (\mathbb{R},0)$ has an isolated stratified critical point at $0$, we can assume that $\Gamma_{g,\rho}^{X^f}$ is a $C^1$-subanalytic curve in a neighborhood of the origin. 
\begin{lemma}\label{Lemma2}
We have $\Gamma_{g,\rho}^{X^f} \cap X^g \subset \{0\}$ in a neighborhood of the origin.
\end{lemma}
\proof Same proof as Lemma \ref{Lemma1}. \endproof

Therefore if $\epsilon >0$ is small enough, $g_{\vert \tilde{S}_\epsilon \cap X^f}$ has a finite number of critical points. They do not lie in the level $\{g=0 \}$ so by Lemma \ref{CrPts2}, they are outwards-pointing for $g_{\vert X^f \cap \tilde{B}_\epsilon}$ if they lie in $\{ g > 0 \}$ and inwards-pointing if they lie in $\{g < 0 \}$. 

Let us study now the critical points of $g_{\vert X \cap f^{-1}(\delta) \cap \tilde{S}_\epsilon }$. We will need the following lemma.
\begin{lemma}\label{Lemma3}
For every $\epsilon >0$ sufficiently small, there exists $\delta(\epsilon)>0$ such that for $0 <\vert  \delta \vert \le \delta(\epsilon)$, the points $p_i^{\delta,\epsilon}$ lie in $\tilde{B}_{\epsilon/4}$.
\end{lemma} 
\proof Let 
$$W = \left\{ (x,r,y) \in U \times \mathbb{R} \times \mathbb{R} \ \vert \ \rho(x)=r, y=f(x) \hbox{ and } x \in \overline{\Gamma_{f,g}} \right\}.$$
Then $W$ is a subanalytic set of  $\mathbb{R}^n \times \mathbb{R} \times \mathbb{R}$ and since it is a graph over $\overline{\Gamma_{f,g}}$, its dimension is less or equal to $1$. Let 
$$\begin{array}{ccccc}
 \pi & : &  \mathbb{R}^n \times \mathbb{R} \times \mathbb{R} &  \rightarrow  & \mathbb{R} \times \mathbb{R} \cr
    &   & (x,r,y) & \mapsto & (r,y), \cr
\end{array}$$
be the projection on the last two factors. Then $\pi_{\vert W} : W \rightarrow \pi(W)$ is proper and $\pi(W)$ is a closed subanalytic set in a neighborhood of the origin.

Let us write $Y_1 = \mathbb{R} \times \{0\}$ and let $Y_2$ be the closure of $\pi(W) \setminus Y_1$. Since $Y_2$ is a curve for $W$ is a curve, $0$ is isolated in $Y_1 \cap Y_2$. By Lojasiewicz's inequality, there exists a constant $C>0$ and an integer $N>0$ such that $ \vert y \vert \ge C r^N$ for $(r,y)$ in $Y_2$ sufficiently close to the origin. So if $x \in \Gamma_{f,g}$ then $\vert f(x) \vert \ge C \rho(x)^N$ if $\rho(x)$ is small enough. 

Let us fix $\epsilon >0$ small. If $0 < \vert \delta \vert \le \frac{1}{C} ( \frac{\epsilon}{4} )^N$ and $x \in f^{-1}( \delta) \cap \Gamma_{f,g}$ then $\rho(x) \le \frac{\epsilon}{4}$. \endproof

For each stratum $V \nsubseteq X^f$, let 
$$\Gamma^V_{f,g,\rho} = \left\{ x \in V \ \vert \ {\rm rank} [ \nabla f_{\vert V}(x), \nabla g_{\vert V}(x), \nabla \rho_{\vert V} (x) ] <3 \right\},$$
and let $\Gamma_{f,g,\rho} = \cup_{V \nsubseteq X^f} \Gamma^V_{f,g,\rho}$. By Lemma \ref{CrPts5}, we can assume that $\Gamma_{f,g,\rho} \setminus \Gamma_{f,g}$ is a $C^1$-subanalytic manifold of dimension 2. 
Let us choose $\epsilon >0$ small enough so that $\tilde{S_\epsilon}$ intersects $\Gamma_{f,g,\rho} \setminus \Gamma_{f,g}$  transversally. Therefore $(\Gamma_{f,g,\rho} \setminus \Gamma_{f,g} ) \cap \tilde{S_\epsilon}$ is subanalytic curve. By Lemma \ref{Lemma2}, we can find $\delta(\epsilon)>0$ such that $f^{-1} \big([\delta(\epsilon),-\delta(\epsilon) ] \big) \cap \tilde{S}_\epsilon \cap \Gamma_{f,g}$ is empty and so 
$$f^{-1} \big([-\delta(\epsilon),\delta(\epsilon) ]\big) \cap \left( \Gamma_{f,g,\rho} \setminus \Gamma_{f,g} \right) \cap \tilde{S}_\epsilon = f^{-1} \big([-\delta(\epsilon),\delta(\epsilon) ] \big) \cap \Gamma_{f,g,\rho}  \cap \tilde{S}_\epsilon .$$
Let $C_1,\ldots,C_t$ be the connected components of $ f^{-1} \big([-\delta(\epsilon),\delta(\epsilon) ] \big) \cap \Gamma_{f,g,\rho}  \cap \tilde{S}_\epsilon$ whose closure intersects $X^f \cap \tilde{S_\epsilon}$. Note that by Thom's $(a_f)$-condition, for each $i \in \{1,\ldots,t \}$, $\overline{C_i} \cap X^f$ is subset of $\Gamma^{X^f}_{g,\rho}$. Let $z_i$ be a point in $\overline{C_i} \cap X^f$. Since $C_i \cap X^f = \emptyset$, there exists $0< \delta_i'(\epsilon) \le \delta(\epsilon)$ such that the fibre $f^{-1}( \delta)$, $0< \vert \delta \vert \le \delta_i'(\epsilon)$, intersects $C_i$ transversally in a neighborhood of $z_i$. 

Let us choose $\delta$ such that $0 <\vert \delta \vert \le {\rm Min} \{ \delta_i'(\epsilon) \ \vert \ i=1,\ldots,t \}$. Then the fibre $f^{-1}( \delta)$ intersect the $C_i$'s transversally and $f^{-1}(\delta) \cap (\cup_i C_i)$ is exactly the set of critical points of $g_{\vert f^{-1}(\delta) \cap X \cap \tilde{S}_\epsilon }$. We have proved:
\begin{lemma}\label{Lemma4}
For $0 < \vert \delta \vert \ll \epsilon \ll 1$, $g_{\vert f^{-1}(\delta) \cap X \cap \tilde{S_\epsilon} }$  has a finite number of critical points, which are exactly the points in $\Gamma_{f,g,\rho} \cap \tilde{S_\epsilon} \cap f^{-1}(\delta) $.
\end{lemma}
$\hfill \Box$

Let $\{s_1^{\delta,\epsilon},\ldots,s_u^{\delta,\epsilon} \}$ be the set of critical points of $g_{\vert f^{-1}(\delta) \cap X \cap \tilde{S_\epsilon}}$.
\begin{lemma}
For $i \in \{1,\ldots,u \}$, $g(s_i^{\delta,\epsilon}) \not= 0$ and $s_i^{\delta,\epsilon}$ is outwards-pointing (resp. inwards-pointing) if and only if $g(s_i^{\delta,\epsilon}) >0$ (resp. $g(s_i^{\delta,\epsilon}) <0$).
\end{lemma}
\proof Note that $s_i^{\delta,\epsilon}$ is necessarily outwards-pointing or inwards-pointing because $s_i^{\delta,\epsilon} \notin \Gamma_{f,g}$. 

Assume that for each $\delta >0$ small enough, there exists a point $s_i^{\delta,\epsilon}$ such that $g(s_i^{\delta,\epsilon})=0$. Then we can construct a sequence of points $(\sigma_n)_{n \in \mathbb{N}}$ such that $g(\sigma_n)=0$ and $\sigma_n$ is a critical point of $g_{\vert f^{-1}(\frac{1}{n}) \cap X \cap \tilde{S}_\epsilon }$. We can also assume that the points $\sigma_n$'s belong to the same stratum $S$ and that they tend to  $\sigma \in V$ where $V \subseteq X^f$ and $V \subset \partial \overline{S}$.  Therefore we have a decomposition:
$$\nabla g_{\vert S}( \sigma_n)= \lambda_n \nabla f_{\vert S} (\sigma_n) + \mu_n \nabla \rho_{\vert S} (\sigma_n).$$
Now by Whitney's condition (a), $T_{\sigma_n} S $ tends to a linear space $T$ such that $T_\sigma V \subset T$. So $\nabla g_{\vert S} (\sigma_n)$ tends to a vector in $T$ whose orthogonal projection on $T_\sigma V$ is exactly $\nabla g_{\vert V} (\sigma)$. Similarly $\nabla \rho_{\vert S} (\sigma_n)$ tends to a vector in $T$ whose orthogonal projection on $T_\sigma V$ is exactly $\nabla \rho_{\vert V} (\sigma)$. By Thom's condition, $\nabla  f_{\vert S} (\sigma_n)$ tends to a vector in $T$ which is orthogonal to $T_\sigma V$, so we see that $\nabla g_{\vert V} (\sigma)$ and $\nabla \rho_{\vert V}(\sigma)$ are colinear which means that $\sigma$ is a critical point of $g_{\vert X^f \cap \tilde{S}_\epsilon}$. But since $g(\sigma_n)=0$, we find that $g(\sigma)=0$, which is impossible by Lemma \ref{Lemma2}. This proves the first assertion.

To prove the second one, we use the same method. Assume that for each $\delta >0$ small enough, there exists a point $s_i^{\delta,\epsilon}$ such that $g(s_i^{\delta,\epsilon}) >0$ and $s_i^{\delta,\epsilon}$ is an inwards-pointing critical point for $g_{\vert X \cap f^{-1}(\delta) \cap \tilde{S}_\epsilon}$.  Then we can construct a sequence of points $(\tau_n)_{n \in \mathbb{N}}$ such that $g(\tau_n) >0$ and $\tau_n$ is an inwards-pointing critical point for $g_{\vert  f^{-1}(\frac{1}{n}) \cap  X \cap \tilde{S}_\epsilon}$. We can also assume that the points $\tau_n$'s belong to the same stratum $S$ and that they tend to  $\tau \in V$ where $V \subseteq X^f$ and $V \subset \partial \overline{S}$. Therefore, we have a decomposition:
$$\nabla g_{\vert S}( \tau_n)= \lambda_n \nabla f_{\vert S} (\tau_n) + \mu_n \nabla \rho_{\vert S} (\tau_n),$$
with $\mu_n <0$. Using the same arguments as above, we find that $\nabla g_{\vert V} (\tau)= \mu \nabla \rho_{\vert S} (\tau)$ with $\mu \le 0$ and $g(\tau) \ge 0$. This contradicts the remark after Lemma \ref{Lemma2}. 
Of course, this proof works for $\delta <0$.
\endproof

Let $\Gamma_{g,\rho}$ be the following polar set:
$$\Gamma_{g,\rho} = \left\{ x \in U \ \vert \ {\rm rank} [ \nabla g (x) , \nabla \rho (x) ] < 2  \right\}.$$
By Lemma \ref{CrPts4} and Lemma \ref{CrPts1}, we can assume that $\Gamma_{g,\rho} \setminus \{g=0\}$ does not intersect $X^f \setminus \{0\}$ in  a neighborhood of $0$ and so $\Gamma_{g,\rho} \setminus \{g=0\}$ does not intersect $X^f \cap \tilde{S}_\epsilon$ for $\epsilon >0$ sufficiently small. Since the critical points of $g_{\vert X^f \cap \tilde{S}_\epsilon}$ lie outside $\{g=0 \}$, they do not belong to $\Gamma_{g,\rho}\cap \tilde{S}_\epsilon$ and so the critical points of $g_{\vert f^{-1}( \delta) \cap X \cap \tilde{S_\epsilon}}$ do not neither if $\delta $ is sufficiently small. Hence at each critical point of $g_{\vert f^{-1}( \delta) \cap X \cap \tilde{S_\epsilon}}$, $g_{\vert \tilde{S_\epsilon}}$ is a submersion. We are in position to apply Theorem 3.1 and Lemma 2.1 in \cite{DutertreManuscripta}. For $0<  \vert \delta \vert \ll \epsilon \ll 1$, we set 
$$I(\delta, \epsilon, g) = \sum_{i=1}^r {\rm ind}(g, X \cap f^{-1}(\delta),p_i ^{\delta,\epsilon}),$$
$$I(\delta, \epsilon, -g) = \sum_{i=1}^r {\rm ind}(-g,X \cap f^{-1}(\delta), p_i ^{\delta,\epsilon}).$$
\begin{theorem}
We have 
$$\displaylines{
\quad I( \delta,\epsilon, g) + I( \delta,\epsilon,-g) = 2 \chi \big(X \cap f^{-1}( \delta) \cap \tilde{B}_\epsilon \big)  \hfill \cr 
\hfill-\chi \big( X \cap f^{-1}(\delta) \cap \tilde{S}_\epsilon \big) -
 \chi \big( X^g \cap f^{-1}( \delta) \cap \tilde{S}_\epsilon \big). \quad \cr
}$$
\end{theorem}
\proof Let us denote by $\{a_j^+\}_{j=1}^{\alpha^+}$ (resp. $\{a_j^-\}_{j=1}^{\alpha^-}$) the outwards-pointing (resp. inwards-pointing) critical points of $ g : X \cap f^{-1}(\delta) \cap \tilde{S_\epsilon} \rightarrow \mathbb{R}$. Applying Morse theory type theorem (\cite{DutertreManuscripta}, Theorem 3.1) and using Lemma 2.1 in \cite{DutertreManuscripta}, we can write
$$I(\delta,\epsilon,g) + \sum_{j=1}^{\alpha^-} {\rm ind}(g,X \cap f^{-1}(\delta) \cap \tilde{S}_\epsilon, a_j^-) = \chi (X \cap f^{-1}(\delta) \cap \tilde{B}_\epsilon ) \eqno(1),$$
$$I(\delta,\epsilon,-g) + \sum_{j=1}^{\alpha^+}{\rm ind}(-g,X \cap f^{-1}(\delta) \cap \tilde{S}_\epsilon,-a_j^+) = \chi (X \cap f^{-1}(\delta) \cap \tilde{B}_\epsilon ) \eqno(2).$$
Let us evaluate 
$$\sum_{j=1}^{\alpha^-} {\rm ind}(g,X \cap f^{-1}(\delta) \cap \tilde{S}_\epsilon, a_j^-) +  \sum_{j=1}^{\alpha^+} {\rm ind}(-g,X \cap f^{-1}(\delta) \cap \tilde{S}_\epsilon,a_j^+).$$
Since the outwards-pointing critical points of $g_{\vert X \cap f^{-1}(\delta) \cap \tilde{S}_\epsilon }$ lie in $\{ g >0 \}$ and the inwards-pointing critical points of $g_{\vert X \cap f^{-1}(\delta) \cap \tilde{S}_\epsilon }$ lie in $\{ g <0 \}$, we have
$$\displaylines{
\quad \chi ( X \cap f^{-1}(\delta) \cap \tilde{S}_\epsilon \cap \{ g \ge 0 \}) - \chi ( X \cap f^{-1}(\delta) \cap \tilde{S}_\epsilon \cap \{ g = 0 \}) = \hfill \cr
\hfill \sum_{j=1}^{\alpha^+} {\rm ind}(g,X \cap f^{-1}(\delta) \cap \tilde{S}_\epsilon,,a_j^+)   \quad  (3),\cr
}$$
and
$$\displaylines{
\quad \chi ( X \cap f^{-1}(\delta) \cap \tilde{S}_\epsilon \cap \{ g \le 0 \}) - \chi ( X \cap f^{-1}(\delta) \cap \tilde{S}_\epsilon \cap \{ g = 0 \}) =  \hfill \cr
\hfill \sum_{j=1}^{\alpha^-} {\rm ind}(-g,X \cap f^{-1}(\delta) \cap \tilde{S}_\epsilon,a_j^-) \quad (4). \cr
}$$
Therefore making $(3)+(4)$ and using the Mayer-Vietoris sequence, we find
$$\displaylines{
\quad \chi( X \cap f^{-1}(\delta) \cap \tilde{S}_\epsilon)- \chi( X \cap f^{-1}(\delta) \cap \tilde{S}_\epsilon \cap \{g=0 \}) = \hfill \cr
\hfill  \sum_{j=1}^{\alpha^+} {\rm ind}(g,X \cap f^{-1}(\delta) \cap \tilde{S}_\epsilon,a_j^+) + \sum_{j=1}^{\alpha^-} {\rm ind}(-g,X\cap f^{-1}(\delta) \cap \tilde{S}_\epsilon,a_j^-) \quad (5). \cr
}$$
Moreover we have
$$\displaylines{
\quad \chi (X \cap f^{-1}(\delta) \cap \tilde{S}_\epsilon ) = \sum_{j=1}^{\alpha^+} {\rm ind}(g,X \cap f^{-1}(\delta) \cap \tilde{S}_\epsilon,a_j^+)  \hfill \cr
\hfill + \sum_{j=1}^{\alpha^-} {\rm ind}(g,X\cap f^{-1}(\delta) \cap \tilde{S}_\epsilon,a_j^-) \quad (6), \cr
}$$
$$\displaylines{
\quad \chi (X \cap f^{-1}(\delta) \cap \tilde{S}_\epsilon ) = \sum_{j=1}^{\alpha^+} {\rm ind}(-g,X \cap f^{-1}(\delta) \cap \tilde{S}_\epsilon,a_j^+)  \hfill \cr
\hfill + \sum_{j=1}^{\alpha^-} {\rm ind}(-g,X\cap f^{-1}(\delta) \cap \tilde{S}_\epsilon,a_j^-) \quad (7). \cr
}$$
The combination $-(5)+(6)+(7)$ leads to
$$\displaylines{
\quad \chi(X \cap f^{-1}(\delta) \cap \tilde{S}_\epsilon ) + \chi(X \cap f^{-1}(\delta) \cap \tilde{S}_\epsilon \cap \{g=0 \}) = \hfill \cr
\hfill \sum_{j=1}^{\alpha^+} {\rm ind}(-g,X \cap f^{-1}(\delta) \cap \tilde{S}_\epsilon,a_j^+) + \sum_{j=1}^{\alpha^-} {\rm ind}(g,X\cap f^{-1}(\delta) \cap \tilde{S}_\epsilon,a_j^-). \quad \cr
}$$
\endproof

Let us assume now that  $(X,0)$ is equipped with a Whitney stratification $\mathcal{W}= \cup_{\alpha \in A}  W_\alpha $ and $f : (X,0) \rightarrow (\mathbb{R},0)$ has an isolated critical point at $0$. In this situation, our results apply taking for $\mathcal{V}$ the following stratification:
$$\left\{ W_\alpha \setminus f^{-1}(0), W_\alpha \cap f^{-1}(0) \setminus \{0\}, \{0\} \ \vert \  W_\alpha \in \mathcal{W} \right\}.$$
\begin{corollary}
If $f : (X,0) \rightarrow (\mathbb{R},0)$ has an isolated stratified critical point at $0$, then 
$$\displaylines{
\quad  I(\delta,\epsilon,g) + I(\delta,\epsilon,-g) = 2 \chi \big(X \cap f^{-1}( \delta) \cap \tilde{B}_\epsilon \big)  \hfill \cr
\hfill- \chi \big(X^f \cap \tilde{S}_\epsilon \big) - \chi \big( X^f \cap X^g \cap \tilde{S}_\epsilon \big). \quad \cr
}$$
\end{corollary}
\proof For each stratum $W$ of $X$, let 
$$\Gamma_{f,\rho}^W = \left\{ x \in W \ \vert \ {\rm rank} [ \nabla f_{\vert W}(x), \nabla \rho_{\vert W}(x) ] < 2 \right\},$$
and let $\Gamma_{f,\rho}= \cup_{W} \Gamma_{f,\rho}^W$. 
By Lemma \ref{Lemma2} applied to $X$ and $f$ instead of $X^f$ and $g$, $\Gamma_{f,\rho} \cap \{f=0\} \subset \{0\}$ in a neighborhood of the origin and so $0$ is a regular value of $f : X \cap \tilde{S}_\epsilon \rightarrow \mathbb{R}$ for $\epsilon$ sufficiently small. 
By Thom-Mather's second isotopy lemma, $X \cap f^{-1}(0) \cap \tilde{S}_\epsilon$ is homeomorphic to $X \cap f^{-1}(\delta) \cap \tilde{S}_\epsilon$ for $\delta$ sufficiently small. 

Now let $p$ be a stratified critical point of $f : X^g \rightarrow \mathbb{R}$. By Lemma \ref{CrPts1}, we know that $p$ belongs to $f^{-1}(0) \cap X^g$ and so $p$ is also a critical point of $g : X^f \rightarrow \mathbb{R}$. Hence $p=0$ by Condition (A), and $f : X^g \rightarrow \mathbb{R}$ has an isolated stratified critical point at $0$. As above, we conclude that $X^f \cap X^g \cap \tilde{S}_\epsilon$ is homeomorphic to $X^g \cap f^{-1}( \delta) \cap \tilde{S}_\epsilon$. \endproof

Let $\omega(x) = \sqrt{x_1^2+\cdots+x_n^2}$ be the euclidian distance to the origin. As explained by Durfee in \cite{Durfee}, Lemma 1.8 and Lemma 3.6, there is a neighborhood $\Omega$ of $0$ in $\mathbb{R}^n$ such that for every stratum $V$ of $X^f$, $\nabla \omega_{\vert V}$ and $\nabla \rho_{\vert V}$ are non-zero and do not point in opposite direction in $\Omega \setminus \{0\}$. Applying Durfee's argument (\cite{Durfee}, Proposition 1.7 and Proposition 3.5), we see that $X^f \cap \tilde{S}_\epsilon$ is homeomorphic to $X^f \cap S_{\epsilon'}$ for $\epsilon,\epsilon'>0$ sufficiently small. Similarly $X^f \cap X^g \cap \tilde{S}_\epsilon$ and 
$X^f \cap X^g \cap S_{\epsilon'}$ are homemorphic. Now let us compare $X \cap f^{-1}(\delta) \cap \tilde{B}_\epsilon$ and $X \cap f^{-1}(\delta) \cap B_{\epsilon'}$. Let us choose $\epsilon'$ and $\epsilon$ such that $$X \cap f^{-1}( \delta) \cap B_{\epsilon'} \subset X \cap f^{-1}( \delta) \cap \tilde{B}_\epsilon \subset \Omega.$$
If $\delta $ is sufficiently small then, for every stratum $V \nsubseteq X^f$,  $\nabla \omega_{\vert V \cap f^{-1}( \delta)}$ and $\nabla \rho_{\vert V \cap f^{-1}( \delta)}$ are non-zero and do not point in opposite direction in $\tilde{B}_\epsilon \setminus \mathring{B_{\epsilon'}}$. Otherwise, by Thom's ($a_f$)-condition, we would find a point $p$ in $X^f \cap (\tilde{B}_\epsilon \setminus \mathring{B_{\epsilon'}})$ such that either $\nabla \omega_{\vert S}(p)$ or $\nabla \rho_{\vert S}(p)$ vanish or $\nabla \omega_{\vert S}(p)$ and $\nabla \rho_{\vert S}(p)$ point in opposite direction, where $S$ is the stratum of $X^f$ that contains $p$. This is impossible if we are sufficiently close to the origin. Now, applying the same arguments as Durfee \cite{Durfee}, Proposition 1.7 and Proposition 3.5, we see that $X \cap f^{-1}( \delta) \cap \tilde{B}_\epsilon$ is homeomorphic to $X \cap f^{-1}( \delta) \cap B_{\epsilon'}$ and that $X \cap f^{-1}( \delta) \cap \tilde{S}_\epsilon$ is homeomorphic to $X \cap f^{-1}(\delta) \cap S_{\epsilon'}$. 
\begin{theorem}\label{MainTh}
We have
$$ I(\delta,\epsilon,g)+I(\delta,\epsilon,-g) = 2 \chi(M_f^{\delta,\epsilon}) 
- \chi (X \cap f^{-1}(\delta) \cap S_\epsilon)- \chi ( X^g \cap f^{-1}(\delta) \cap S_\epsilon).$$
\end{theorem}
$\hfill \Box$

\begin{corollary}\label{MainCor}
If $f : (X,0) \rightarrow (\mathbb{R},0)$ has an isolated stratified critical point at $0$, then 
$$I(\delta,\epsilon,g) + I(\delta,\epsilon,-g) = 2 \chi ( M_f^{\delta,\epsilon}) - \chi ({\rm Lk}(X^f)) - \chi ({\rm Lk} (X^f \cap X^g)).$$
\end{corollary}
$\hfill \Box$

Let us remark if ${\rm dim}\ X=2$ then in Theorem \ref{MainTh} and in Corollary \ref{MainCor}, the last term of the right-hand side of the equality vanishes. If ${\rm dim} \ X=1$ then in Theorem \ref{MainTh} and in Corollary \ref{MainCor}, the last two terms of the right-hand side of the equality vanish. 

%Let us apply these results when $X = \mathbb{R}^n$. In this case, $f^{-1}(\delta) \cap B_\epsilon$ is a manifold with boundary (possibly empty) of dimension $n-1$, $f^{-1}(\delta) \cap S_\epsilon$ is a compact manifold of dimension $n-2$ and, since $g : (\mathbb{R}^n,0) \rightarrow (\mathbb{R},0)$ has an isolated critical point at the origin, $g^{-1}(0) \cap f^{-1}(\delta) \cap S_\epsilon$ is a compact manifold of dimension $n-3$, boundary of the manifold $g^{-1}(0) \cap f^{-1}(\delta) \cap B_\epsilon$. Furthermore, if $p$ is a critical point of $g_{\vert f^{-1}(\delta) \cap B_\epsilon}$ then $${\rm ind}(f^{-1}(\delta),g,p)=(-1)^{n-1} {\rm ind}(f^{-1}(\delta),-g,p).$$ Therefore if $n$ is odd, then $$I(\delta,\epsilon,g)+I(\delta,\epsilon,-g) = 2 I(\delta,\epsilon,g),$$ $$\chi (X \cap f^{-1}(\delta) \cap S_\epsilon)=0,$$ and $$\chi(X^g \cap f^{-1}(\delta) \cap S_\epsilon)=2 \chi( X^g \cap f^{-1}(\delta) \cap B_\epsilon).$$ Theorem \ref{MainTh} becomes $$I(\delta,\epsilon,g)= \chi(f^{-1}(\delta) \cap B_\epsilon)- \chi ( f^{-1}(\delta) \cap \{f=0 \} \cap B_\epsilon).$$
%If $n$ is even, then $I(\delta,\epsilon,g) + I(\delta,\epsilon,-g)=0$, 
%$$2 \chi ( f^{-1}(\delta) \cap B_\epsilon)-\chi (f^{-1}(\delta) \cap S_\epsilon)=0,$$ and $\chi (f^{-1}(\delta) \cap \{ g=0 \} \cap S_\epsilon)=0$, so in this setting, our formula is trivially true. 

\section{An infinitesimal Gauss-Bonnet formula}
In this section, we apply the results of the previous section to the case of linear forms and we establish a Gauss-Bonnet type formula for the real Milnor fibre.

We will first show that generic linear forms satisfy Condition (A) and Condition (B). For $v \in S^{n-1}$, let us denote by $v^*$ the function $v^*(x)= \langle v,x \rangle$.
\begin{lemma}\label{LIN1}
There exists a subanalytic set $\Sigma_1 \subset S^{n-1}$ of positive codimension such  that if $v \notin \Sigma_1$, $\{v^*=0 \}$ intersects $X\setminus \{0\}$ transversally (in the stratified sense) in a neighborhood of the origin.
\end{lemma}
\proof It is a particular case of Lemma 3.8 in \cite{DutertreGeoDedicata}. \endproof

\begin{corollary}\label{LIN2}
If $v \notin \Sigma_1$ then $v^*_{\vert X} : (X,0) \rightarrow (\mathbb{R},0)$ has an isolated stratified point at $0$.
\end{corollary}
\proof By Lemma \ref{CrPts1}, we know that the stratified critical points of $v^*_{\vert X}$ lie in $\{v^*=0\}$. But since $\{v^*=0\}$ intersects $X \setminus \{0\}$ transversally, the only possible critical point of $v^*_{\vert X} : (X,0) \rightarrow (\mathbb{R},0)$ is the origin. \endproof

\begin{lemma}\label{LIN3}
There exists a subanalytic set $\Sigma_2 \subset S^{n-1}$ of positive codimension such that if $v \notin \Sigma_2$, then $\Gamma_{f,v^*}$ is a $C^1$-subanalytic curve (possibly empty) in a neighborhood of $0$.
\end{lemma}
\proof Let $V$ be stratum of dimension $e$ such that $V \nsubseteq X^f$. We can assume that $e \ge 2$. Let 
$$M_V = \left\{ (x,y) \in V  \times \mathbb{R}^n \ \vert \  {\rm rank} [ \nabla f_{\vert V} (x), \nabla  y^*_{\vert V}(x) ] < 2 \right\}.$$
It is a subanalytic manifold of class $C^1$ and of dimension $n+1$. To see this, let us pick a point $(x,y)$ in $ M_V$. In a neighborhood of $x$, $V$ is defined by the vanishing of $k=n-e$ $C^2$-subanalytic functions $f_1,\ldots,f_k$. Since $V$ is not included in $X^f$, $f : V \rightarrow \mathbb{R}$ is a submersion and we can assume that in a neighborhood of $x$, the following $(k+1) \times (k+1)$-minor:
$$\frac{\partial(f_1,\ldots,f_k,f)}{\partial(x_1,\ldots,x_k,x_{k+1})},$$
does not vanish. Therefore, in a neighborhood of $(x,y)$, $M_V$ is defined by the vanishing of the following $(k+2) \times (k+2)$-minors:
$$\frac{\partial(f_1,\ldots,f_k,f,y^*)}{\partial(x_1,\ldots,x_k,x_{k+1},x_{k+2})},\ldots,\frac{\partial(f_1,\ldots,f_k,f,y^*)}{\partial(x_1,\ldots,x_k,x_{k+1},x_{n})}.$$
A simple computation of determinants shows that the gradient vectors of these minors are linearly independent. As in  previous lemmas, we show that $\Sigma_{f,v^*}$ is one-dimensional considering the projection
$$\begin{array}{ccccc}
\pi_2 & :  & M^V & \rightarrow & \mathbb{R}^n \cr
  &  & (x,y) & \mapsto & y. \cr
\end{array}$$
Since $\Gamma_{f,v^*} = \cup_{V \nsubseteq X^f} \Gamma_{f,v^*}^V $, we get the result. \endproof
Let $\Sigma= \Sigma_1 \cup \Sigma_2$, it is a subanalytic subset of $S^{n-1}$ of positive codimension and if $v \notin \Sigma$ then $v^*$ satisfies Conditions (A) and (B). In particular, $v^*_{\vert f^{-1}(\delta) \cap X \cap B_\epsilon}$  has a finite number of critical points $p_1^{\delta,\epsilon}, \ldots,p_{r_v}^{\delta,\epsilon}$.
We recall that
$$I(\delta,\epsilon,v^*)= \sum_{i=1}^{r_v} {\rm ind}(v^*,X \cap f^{-1}(\delta),p_i^{\delta,\epsilon}),$$
$$I(\delta,\epsilon,-v^*)= \sum_{i=1}^{r_v} {\rm ind}(-v^*,X \cap f^{-1}(\delta), p_i^{\delta,\epsilon}).$$
In this situation, Theorem \ref{MainTh} and Corollary \ref{MainCor} become
\begin{corollary}\label{MainThLIN}
If $v \notin \Sigma$ then 
$$I(\delta,\epsilon,v^*)+ I( \delta,\epsilon,-v^*) = 2 \chi (M_f^{\delta,\epsilon})
 - \chi (X \cap f^{-1}( \delta) \cap S_\epsilon) - \chi(X^{v^*} \cap f^{-1}( \delta) \cap S_\epsilon). $$
Furthermore, if $f : (X,0) \rightarrow (\mathbb{R},0)$ has an isolated stratified critical point at $0$, then 
$$I(\delta,\epsilon,v^*) + I(\delta,\epsilon,-v^*) = 2 \chi ( M_f^{\delta,\epsilon}) - \chi ({\rm Lk}(X^f)) - \chi ({\rm Lk} (X^f \cap X^{v^*})).$$
\end{corollary}
$\hfill \Box$

As an application, we give a Gauss-Bonnet formula for the Milnor fibre $M_f^{\delta,\epsilon}$. Let $\Lambda_0(X \cap f^{-1}(\delta), -)$ be the Gauss-Bonnet measure on $X \cap f^{-1}( \delta)$ defined by 
$$\Lambda_0(X \cap f^{-1}(\delta), U') = \frac{1}{s_{n-1}} \int_{S^{n-1}} \sum_{x \in U'} {\rm ind}( v^*,X \cap f^{-1}(\delta), x) dx,$$
where $U'$ is a Borel set of $X \cap f^{-1}(\delta)$ (see \cite {BroeckerKuppe}, page 299) and $s_{n-1}$ is the volume of the unit sphere $S^{n-1}$. Note that if $x$ is not a critical point of $v^*_{\vert X \cap f^{-1}(\delta)}$ then ${\rm ind}(v^*,X \cap f^{-1}(\delta),x )=0$. We are going to evaluate $$\lim_{\epsilon \rightarrow 0} \lim_{\delta \rightarrow 0} \Lambda_0 (X \cap f^{-1}(\delta), M_f^{\delta,\epsilon}).$$
\begin{theorem}\label{GB}
We have
$$\displaylines{
\quad \lim_{\epsilon \rightarrow 0} \lim_{\delta \rightarrow 0} \Lambda_0 (X \cap f^{-1}(\delta), M_f^{\delta,\epsilon}) =  
 \chi (M_f^{\delta,\epsilon}) -\frac{1}{2} \chi (X \cap f^{-1}(\delta) \cap S_\epsilon) \hfill \cr
\hfill  -\frac{1}{2 s_{n-1}} \int_{S^{n-1}} \chi (X \cap f^{-1}(\delta) \cap \{v^*=0 \} \cap S_\epsilon) dv. \quad \cr
}$$
Furthermore, if $f : (X,0) \rightarrow (\mathbb{R},0)$ has an isolated stratified critical point at $0$, then 
$$\displaylines{
\quad \lim_{\epsilon \rightarrow 0} \lim_{\delta \rightarrow 0} \Lambda_0 (X \cap f^{-1}(\delta), M_f^{\delta,\epsilon}) =  
 \chi (M_f^{\delta,\epsilon}) - \frac{1}{2} \chi ({\rm Lk}(X^f))\hfill \cr
\hfill  -\frac{1}{2 s_{n-1}} \int_{S^{n-1}} \chi ({\rm Lk} (X^f \cap X^{v^*})) dv. \quad \cr
}$$
\end{theorem}
\proof By definition, we have 
$$
 \Lambda_0 (X \cap f^{-1}(\delta), M_f^{\delta,\epsilon}) = 
 \frac{1}{s_{n-1}} \int_{S^{n-1}} \sum_{x \in M_f^{\delta,\epsilon}} {\rm ind}(v^*,X \cap f^{-1}(\delta),x) dv. 
$$
It is not difficult to see that 
$$\displaylines{
\quad \Lambda_0 (X \cap f^{-1}(\delta), M_f^{\delta,\epsilon}) = \hfill \cr 
\hfill \frac{1}{2 s_{n-1}} \int_{S^{n-1}} \Big[  \sum_{x \in M_f^{\delta,\epsilon}} {\rm ind}(v^*,X \cap f^{-1}(\delta),x) +
 {\rm ind}(-v^*,X \cap f^{-1}(\delta),x)  \Big]dv. \quad \cr
}$$
Note that if $v \notin \Sigma$ then 
$$\sum_{x \in M_f^{\delta,\epsilon}} {\rm ind}(v^*,X \cap f^{-1}(\delta),x) + {\rm ind}(-v^*,X \cap f^{-1}(\delta),x)$$
is equal to $I(\delta,\epsilon,v^*)+ I (\delta,\epsilon,-v^*)$ and is uniformly bounded by Hardt's theorem. By Lebesgue's theorem, we obtain
$$\displaylines{
\quad \lim_{\epsilon \rightarrow 0} \lim_{\delta \rightarrow 0}  \Lambda_0 (X \cap f^{-1}(\delta), M_f^{\delta,\epsilon}) = \hfill \cr
\hfill  \frac{1}{2 s_{n-1}} \int_{S^{n-1}} 
\lim_{\epsilon \rightarrow 0} \lim_{\delta \rightarrow 0} [I(\delta,\epsilon,v^*) + I(\delta,\epsilon,-v^*) ] dv. \quad \cr
}$$
We just have to apply the previous corollary to conclude. \endproof

\section{Infinitesimal linear kinematic formulas}
In this section, we apply the results of the previous section to the case of a linear function in order to obtain ``infinitesimal" linear kinematic formulas for closed subanalytic germs.

We start recalling known facts on the geometry of subanalytic sets. We need some notations:
\begin{itemize}
\item for $k \in \{0,\ldots,n \}$, $G_{n}^k$ is the Grassmann manifold of $k$-dimensional linear
subspaces in $\mathbb{R}^{n}$ and $g_n^k$ is its volume,
\item for $k \in \mathbb{N}$, $b_k$ is the volume of the $k$-dimensional unit ball and $s_k$ is the
volume of the $k$-dimensional unit sphere.
\end{itemize}
In \cite{Fu}, Fu developed integral geometry for compact subanalytic sets. Using the technology of the normal cycle, he associated with every compact subanalytic set $X \subset \mathbb{R}^n$ a sequence of curvature measures 
$$\Lambda_0(X,-),\ldots,\Lambda_n(X,-),$$
called the Lipschitz-Killing measures. He proved several integral geometry formulas, among them a Gauss-Bonnet formula and a kinematic formula. Later another description of the measures using stratified Morse theory was given by Broecker and Kuppe \cite {BroeckerKuppe} (see also \cite{BernigBroecker}). The reader can refer to \cite{DutertreGeoDedicata}, Section 2, for a rather complete presentation of these two approaches and for the definition of the Lipschitz-Killing measures. 

Let us give some comments on these Lipschitz-Killing curvatures. If dim $X=d$ then
$$\Lambda_{d+1}(X,U')=\cdots=\Lambda_n(X,U')=0,$$
for any Borel set $U'$ of $X$ and $\Lambda_d(X,U')= \mathcal{L}_d(U')$, where $\mathcal{L}_d$ is the $d$-dimensional Lebesgue measure in $\mathbb{R}^n$. Furthemore if $X$ is smooth then for any Borel set $U'$ of $X$
and for $k \in \{0,\ldots,d \}$, $\Lambda_k (X,U')$ is related to the classical Lipschitz-Killing-Weil curvature $K_{d-k}$ through the following equality:
$$\Lambda_k (X,U')= \frac{1}{s_{n-d-k-1}} \int_{U'} K_{d-k}(x) dx.$$
In \cite{DutertreGeoDedicata}, Section 5, we studied the asymptotic behaviour of the Lipschitz-Killing measures in the neighborhood of a point of $X$. Namely we proved the following theorem (\cite{DutertreGeoDedicata}, Theorem 5.1).
\begin{theorem}\label{CurvAndLink}
Let $X \subset \mathbb{R}^n$ be a closed subanalytic set such that $0 \in X$. 
We have:
$$\lim_{\epsilon \rightarrow 0} \Lambda_0(X,X\cap B_\epsilon)=1-\frac{1}{2} \chi (\hbox{\em Lk}(X))-\frac{1}{2g_n^{n-1}} \int_{G_n^{n-1}} \chi (\hbox{\em Lk}(X \cap H)) dH.$$
Furthermore for $k \in \{1,\ldots,n-2 \}$, we have:
$$\displaylines{
\qquad \lim_{\epsilon \rightarrow 0} \frac{\Lambda_k(X,X\cap B_\epsilon)}{b_k \epsilon^k} =- \frac{1}{2g_n^{n-k-1}} \int_{G_n^{n-k-1}} \chi(\hbox{\em Lk}(X \cap H)) dH \hfill \cr
\hfill + \frac{1}{2g_n^{n-k+1}} \int_{G_n^{n-k+1}} \chi(\hbox{\em Lk}(X \cap L)) dL, \qquad  \cr
}$$
and:
$$\lim_{\epsilon \rightarrow 0} \frac{\Lambda_{n-1}(X,X\cap B_\epsilon)}{b_{n-1} \epsilon^{n-1}} = \frac{1}{2g_n^2} \int_{G_n^2} \chi(\hbox{\em Lk}(X \cap H)) dH,$$
$$\lim_{\epsilon \rightarrow 0} \frac{\Lambda_{n}(X,X\cap B_\epsilon)}{b_{n} \epsilon^{n}} = \frac{1}{2g_n^1} \int_{G_n^1} \chi(\hbox{\em Lk}(X \cap H)) dH.$$
\end{theorem}

In the sequel, we will use these equalities and Theorem \ref{GB} to establish linear kinematic types formulas for the quantities $\lim_{\epsilon \rightarrow 0} \frac{\Lambda_k(X,X \cap B_\epsilon)}{b_k \epsilon^k}$, $k=1,\ldots,n$. Let us start with some lemmas. We work with a closed subanalytic set $X$ such that $0 \in X$, equipped with a Whitney stratification $\{W_\alpha\}_{\alpha \in A}$. 
\begin{lemma}\label{LemmaLink} 
Let $f$ be a $C^2$-subanalytic function such that $f_{\vert X} : X \rightarrow \mathbb{R}$ has an isolated stratified critical point at $0$. Then for $0< \delta \ll \epsilon \ll 1$, we have
$$\chi (M_f^{\delta,\epsilon})  + \chi (M_f^{-\delta,\epsilon} )= \chi ({\rm Lk}(X)) + \chi ({\rm Lk}(X^f)) .$$
\end{lemma}
\proof With the same technics and arguments as the ones we used in order to establish Corollary \ref{MainCor}, we can prove that
$${\rm ind}(f,X,0) + {\rm ind}(-f,X,0) = 2 \chi (X \cap B_\epsilon) - \chi({\rm Lk}(X)) -\chi({\rm Lk}(X^f)).$$
We conclude thanks to the following equalities
$${\rm ind}(f,X,0) = 1 -\chi( M_f^{-\delta,\epsilon}) , \  {\rm ind}(-f,X,0) = 1 -\chi(M_f^{\delta,\epsilon}),$$
and $$\chi (X \cap B_\epsilon)=1.$$
\endproof

\begin{corollary}
There exist a subanalytic set $\Sigma_1 \subset S^{n-1}$ of positive codimension such that if $v \notin \Sigma$ then for $0< \delta \ll \epsilon \ll 1$, 
$$ \chi ( M_{v^*}^{\delta,\epsilon}) + \chi ( M_{v^*}^{-\delta,\epsilon})=  
\chi ({\rm Lk}(X)) + \chi({\rm Lk} (X \cap \{v^*=0 \})). 
$$
\end{corollary}
\proof Apply Corollary \ref{LIN2} and Lemma \ref{LemmaLink}. \endproof

\begin{lemma}\label{Sections}
Let $S \subset \mathbb{R}^n$ be $C^2$-subanalytic manifold. Let $H \in G_n^{n-k}$, $k \in \{1,\ldots,n\}$ and let $G^1_{H^\perp}$ be the Grassmann manifold of lines in the orthogonal complement $H^\perp$ of $H$. There exists a subanalytic set $\Sigma'_H \subset G^1_{H^\perp}$ of positive codimension such that if $\nu \notin \Sigma'_H$ then $H \oplus \nu$ intersects $S \setminus \{0\}$ transversally.
\end{lemma}
\proof Assume that $S$ has dimension $e$ and that $H$ is given by the equations $x_1=\ldots=x_k=0$ so that $H^\perp =\mathbb{R}^k$ with coordinate system $(x_1,\ldots,x_k)$. Let $W$ be defined by
$$\displaylines{
\quad W = \Big\{ (x,v_1,\ldots,v_{k-1}) \in \mathbb{R}^n \times (\mathbb{R}^k )^{k-1} \ \vert \ x \in S\setminus \{0\}  \hfill \cr
\hfill \hbox{ and } \langle x, v_1 \rangle = \cdots = \langle x, v_{k-1} \rangle =0 \Big\},  \quad \cr
}$$
where $v_i \in \mathbb{R}^k \times \{0\} \subset \mathbb{R}^n$. 
Let us show that $W$ is a $C^2$-subanalytic manifold of dimension $e+(k-1)^2$.  Let $(y,w)$ be a point in $W$. We can assume that around $y$, $S$ is defined by the vanishing of $n-e$ $C^2$-subanalytic functions $f_1,\ldots,f_{n-e}$. Hence in a neighborhood of $(y,w)$, $W$ is defined by the equations:
$$f_1(x)=\ldots=f_{n-e}(x)=0 \hbox{ and } 
\langle x,v_1 \rangle = \cdots = \langle x, v_{k-1} \rangle=0.$$ 
Because $y\not=0$, we see that the gradient vectors of this $n-e+k-1$ functions are linearly independent at $(y,w)$. This enables us to conclude that $W$ is a $C^2$-subanalytic manifold of dimension $e+(k-1)^2$. 
Let $\pi_2$ be the following projection:
$$\pi_2 : W \rightarrow (\mathbb{R}^n)^{n-k}, (x,v_1,\ldots,v_{n-k}) \mapsto (v_1,\ldots,v_{n-k}).$$
 Bertini-Sard's theorem implies that the set of critical values of $\pi_2$ is a subanalytic set of positive codimension. 
 If $(v_1,\ldots,v_{k-1})$ lies outside this subanalytic set then the $(n-k+1)$-plane $\{x \in \mathbb{R}^n \ \vert \ \langle x,v_1 \rangle = \cdots = \langle x, v_{k-1} \rangle=0 \}$ contains $H$ and intersects 
 $S \setminus \{0\}$ transversally. $\hfill \Box$

Now we can present our infinitesimal linear kinematic formulas. Let $H \in G_n^{n-k}$, $k \in \{1,\ldots, n \}$, and let $S_{H^\perp}^{k-1}$ be the unit sphere of the orthogonal complement of $H$. Let $v$ be an element in $S_{H^\perp}^{k-1}$. For $\delta >0$, we denote by $H_{v,\delta}$ the $(n-k)$-dimensional affine space $H+\delta v$ and we set 
$$\beta_0(H,v) = \lim_{\epsilon \rightarrow 0} \lim_{\delta \rightarrow 0} \Lambda_0(H_{\delta,v} \cap X, H_{\delta,v} \cap X \cap B_\epsilon).$$
Then we set 
$$\beta_0 (H) =\frac{1}{s_{k-1}}\int_{S_{H^\perp}^{k-1}} \beta_0 (H,v) dv.$$
\begin{theorem}\label{KinForm}
For $k \in \{1,\ldots,n \}$, we have 
$$ \lim_{\epsilon \rightarrow 0} \frac{\Lambda_k(X,X \cap B_\epsilon)}{b_k \epsilon^k } = \frac{1}{g_n^{n-k}} \int_{G_n^{n-k}} \beta_0 (H) dH.$$
\end{theorem}
\proof We treat first the case $k \in \{1,\ldots,n-2\}$. By Theorem \ref{CurvAndLink} , we know that 
$$\displaylines{
\qquad \lim_{\epsilon \rightarrow 0} \frac{\Lambda_k(X,X\cap B_\epsilon)}{b_k \epsilon^k} =- \frac{1}{2g_n^{n-k-1}} \int_{G_n^{n-k-1}} \chi({\rm Lk}(X \cap H)) dH \hfill \cr
\hfill + \frac{1}{2g_n^{n-k+1}} \int_{G_n^{n-k+1}} \chi({\rm Lk}(X \cap L)) dL. \qquad  \cr
}$$
By Lemma 3.8 in \cite{DutertreGeoDedicata}, we know that generically $H$ intersects $X \setminus \{0\}$ transversally in a neighborhood of the origin. Let us fix $H$ that satisfies this generic property. For any $v \in S_{H^\perp}^{k-1}$, let $\nu$ be the line generated by $v$ and let $L_v$ be the $(n-k+1)$-plane defined by $L_v = H \oplus \nu$. By Lemma \ref{Sections}, we know that for $v$ generic in $S_{H^\perp}^{k-1}$, $L_v$ intersects $X \setminus \{0\}$ transversally in a neighborhood of the origin. Therefore, $v^*_{\vert X \cap L_v}$ has an isolated singular point at $0$ and we can apply Theorem \ref{GB}. We have 
$$\displaylines{
\quad \lim_{\epsilon \rightarrow 0} \lim_{\delta \rightarrow 0} \Lambda_0(X \cap L_v \cap \{v^*= \delta \},X \cap L_v \cap \{v^*= \delta \} \cap B_\epsilon ) = \hfill \cr
\quad \quad \chi(X \cap L_v \cap \{v^*=\delta \} \cap B_\epsilon) -\frac{1}{2} \chi ( {\rm Lk}(X \cap L_v \cap \{v^*=0\} ))  \hfill \cr
\hfill -\frac{1}{2 s_{n-k}} \int_{S_{L_v}^{n-k}} \chi ({\rm Lk} (X \cap L_v \cap \{ v^*=0 \} \cap \{w^*=0 \}) )dw, \quad \cr
}$$
where $S_{L_v}^{n-k}$ is the unit sphere of $L_v$.
Let us remark that $L_v \cap \{v^* = \delta \}$ is exactly $H_{v,\delta}$ and that $L_v \cap \{v^*=0\}$ is $H$. We can also apply Lemma \ref{LemmaLink} to $v^*_{\vert X \cap L_v}$ to obtain the following relation:
$$\displaylines{
\quad \beta_0 (H,v)+ \beta_0 (H,-v)  = \chi ({\rm Lk}(X \cap L_v)) \hfill \cr
\hfill -\frac{1}{s_{n-k} }\int_{S_{L_v}^{n-k}} \chi({\rm Lk} (X \cap H \cap  \{w^*=0\} )) dw. \quad \cr 
}$$
Since $\beta(H)$ is equal to 
$$\frac{1}{2 s_{k-1}} \int_{S_{H^\perp}^{k-1} }\left[ \beta_0(H,v)+ \beta_0(H,-v) \right] dv,$$
we find that 
$$\displaylines{
\quad \beta(H)= \frac{1}{2 s_{k-1}} \int_{S_{H^\perp}^{k-1}} \chi ({\rm Lk}(X \cap L_v)) dv  \hfill \cr
\hfill - \frac{1}{2 s_{k-1} s_{n-k}} \int_{S_{H^\perp}^{k-1}} \int_{S_{L_v}^{n-k}} \chi ({\rm Lk}(X \cap H \cap \{ w^*= 0\} ))dw dv. \quad \cr
}$$
Replacing spheres with Grassman manifolds in this equality, we obtain
$$\displaylines{
\quad \beta (H) = \frac{1}{2g_k^1} \int_{G_{H^\perp}^1} \chi ({\rm Lk}(X \cap H \oplus \nu) )d\nu  \hfill \cr
\hfil -\frac{1}{2 g_k^1 g_{n-k+1}^{n-k}} \int_{G_{H^\perp}^1} \int_{G_{H \oplus \nu}^{n-k}} \chi ({\rm Lk} (X \cap H \cap K)) dK d\nu. \quad \cr
}$$ Therefore, we have 
$$\displaylines{
\quad \frac{1}{g_n^{n-k}} \int_{G_n^{n-k}} \beta(H) dH = \frac{1}{2 g_k^1 g_n^{n-k}} \int_{G_n^{n-k}} \int_{G_{H^\perp}^1} \chi ({\rm Lk} (X \cap H \oplus \nu) ) d\nu dH - \hfill \cr
\hfill \frac{1}{2 g_n^{n-k} g_k^1 g_{n-k+1}^{n-k}} 
\int_{G_n^{n-k}} \int_{G_{H^\perp}^1} \int_{G_{H \oplus \nu}^{n-k} }\chi ({\rm Lk} (X \cap H \cap K) ) dK d\nu dH. \quad \cr
}$$
Let us compute 
$$\mathcal{I}= \frac{1}{2g_n^{n-k} g_k^1} \int_{G_n^{n-k}} \int_{G_{H^\perp}^1} \chi ({\rm Lk} (X \cap H \oplus \nu) ) d\nu dH.$$
Let $\mathcal{H}$ be the flag variety of pairs $(L,H)$, $L \in G_n^{n-k+1}$ and $H \in G_L^{n-k}$. This variety is a bundle over $G_n^{n-k}$, each fibre being a $G_k^1$. Hence we have
$$\displaylines{
\quad \int_{G_n^{n-k}} \int_{G_{H^\perp}^1} \chi ({\rm Lk} (X \cap H \oplus \nu)) d\nu dH = \int_{G_n^{n-k+1}} \int_{G_L^{n-k}} \chi ({\rm Lk}(X \cap L)) dH dL = \hfill \cr
\hfill g_{n-k+1}^{n-k} \int_{G_n^{n-k+1}} \chi ({\rm Lk} (X \cap L) )dL. \quad \cr
}$$
Finally, we get that 
$$\displaylines{
\quad \mathcal{I} = \frac{g_{n-k+1}^{n-k}}{2 g_n^{n-k} g_k^1} \int_{G_n^{n-k+1}} \chi ({\rm Lk} (X \cap L)) dL = \hfill \cr
\hfill \frac{1}{2 g_n^{n-k+1}} \int_{G_n^{n-k+1}} \chi ({\rm Lk} (X \cap L)) dL. \quad \cr
}$$
Let us compute now 
$$\mathcal{J} =\frac{1}{2 g_n^{n-k} g_k^1 g_{n-k+1}^{n-k}} \int_{G_n^{n-k}} \int_{G_{H^\perp}^1}  \int_{G_{H \oplus \nu}^{n-k}} \chi({\rm Lk}(X \cap H \cap K)) dK d\nu dH.$$
First, as we have just done above, we can write 
$$\mathcal{J} = \frac{1}{2 g_n^{n-k} g_k^1 g_{n-k+1}^{n-k}} \int_{G_n^{n-k+1}} \int_{G_L^{n-k}} \int_{G_L^{n-k}} \chi( {\rm Lk} (X \cap H \cap K)) dK dH dL.$$
Then we remark (see \cite{DutertreGeoDedicata}, Corollary 3.11 for a similar argument) that 
$$\frac{1}{g_{n-k+1}^{n-k}} \int_{G_L^{n-k}} \chi ({\rm Lk}( X \cap H \cap K)) dK= \frac{1}{g_{n-k}^{n-k-1}} \int_{G_H^{n-k-1}} \chi ({\rm Lk}(X \cap J)) dJ,$$
and so
$$\mathcal{J} = \frac{1}{2g_n^{n-k} g_k^1 g_{n-k}^{n-k-1}} \int_{G_n^{n-k+1}} \int_{G_L^{n-k}} \int_{G_H^{n-k-1}} \chi ({\rm Lk}(X \cap J))  dJ dH dL.$$
Considering the flag variety of pairs $(H,J)$, $H \in G_L^{n-k}$ and $J \in G_H^{n-k-1}$, and proceeding as above, we find
$$\int_{G_L^{n-k}} \int_{G_H^{n-k-1}} \chi ({\rm Lk}(X \cap J)) dJ dH = g_2^1 \int_{G_L^{n-k-1}} \chi ({\rm Lk} (X \cap J)) dJ,$$
so
$$ \mathcal{J}= \frac{g_2^1}{2g_n^{n-k} g_k^1 g_{n-k}^{n-k-1}} \int_{G_n^{n-k+1}} \int_{G_L^{n-k-1}} \chi ({\rm Lk} (X \cap J) )dJ.$$
To finish the computation, we consider the flag variety of pairs $(L,J)$, $L \in G_n^{n-k+1}$ and $J \in G_L^{n-k-1}$. It is a bundle over $G_n^{n-k-1}$, each fibre being a $G_{k+1}^2$. Hence we have
$$\mathcal{J}= \frac{g_2^1}{2 g_n^{n-k} g_k^1 g_{n-k}^{n-k-1}} \int_{G_n^{n-k-1}} \int_{G_{J^\perp}^2} \chi ({\rm Lk}(X \cap J) )dJ d M,$$
$$\displaylines{
\quad \mathcal{J}= \frac{g_2^1g_{k+1}^2}{2 g_n^{n-k} g_k^1 g_{n-k}^{n-k-1}} \int_{G_n^{n-k-1}}  \chi ({\rm Lk}(X \cap J) )dJ = \hfill \cr
\hfill \frac{1}{2 g_n^{n-k-1}} \int_{G_n^{n-k-1}}  \chi ({\rm Lk}(X \cap J)) dJ. \quad \cr
}$$
This ends the proof for the case $k \in \{1,\ldots,n-2 \}$. For $k=n-1$ or $n$, the proof is the same. We just have to remark that in these cases 
$$\beta_0 (H,v) + \beta_0 (H,-v) = \chi ({\rm Lk} (X \cap L_v)),$$
and if $k=n-1$, dim$\ L_v =2$ and if $k=n$, dim$\ L_v=1$. \endproof

Let us end with some remarks on the limits $\lim_{\epsilon \rightarrow 0} \frac{\Lambda_k(X , X \cap B_\epsilon)}{b_k \epsilon^k}$. We already know that if dim$\ X=d$ then $\lim_{\epsilon \rightarrow 0} \frac{\Lambda_k(X , X \cap B_\epsilon)}{b_k \epsilon^k}=0$ for $k \ge d+1$. This is also the case if $l < d_0$, where $d_0$ is the dimension of the stratum that contains $0$. To see this let us first relate the limits $\lim_{\epsilon \rightarrow 0} \frac{\Lambda_k(X , X \cap B_\epsilon)}{b_k \epsilon^k}$ to the polar invariants defined by Comte and Merle in \cite{ComteMerle}. They can be defined as follows. Let $H \in G_n^{n-k}$, $k \in \{1,\ldots, n \}$, and let $v$ be an element in $S_{H^\perp}^{k-1}$. For $\delta >0$, we set 
$$\lambda_0(H,v) = \lim_{\epsilon \rightarrow 0} \lim_{\delta \rightarrow 0} \chi(H_{\delta,v} \cap X \cap B_\epsilon),$$
and then 
$$\sigma_k (X,0)=\frac{1}{s_{k-1} }\int_{S_{H^\perp}^{k-1}} \lambda_0 (H,v) dv.$$
Moreover, we put $\sigma_0 (X,0)=1$. 
\begin{theorem}
For $k \in \{0,\ldots,n-1\}$, we have
$$\lim_{\epsilon \rightarrow 0} \frac{\Lambda_k(X,X \cap B_\epsilon)}{b_k \epsilon^k} = \sigma_k(X,0) -\sigma_{k+1}(X,0).$$
Furthermore, we have
$$\lim_{\epsilon \rightarrow 0} \frac{\Lambda_n(X,X \cap B_\epsilon)}{b_n \epsilon^n} = \sigma_n(X,0).$$
\end{theorem}
\proof It is the same proof as Theorem \ref{KinForm}. For example if $k \in \{0,\ldots,n-1\}$, we just have to remark that
$$\lambda_0 (H,v)+ \lambda_0(H,-v) = \chi ({\rm Lk}(X \cap L_v)) + \chi ({\rm Lk}(X \cap H)),$$
by Lemma \ref{LemmaLink}, which implies that
$$\sigma_k(X,0) = \frac{1}{2 g_n^{n-k+1}} \int_{G_n^{n-k+1}} \chi ({\rm Lk}(X \cap L)) dL + \frac{1}{2 g_n^{n-k}}  \int_{G_n^{n-k}} \chi ({\rm Lk} (X \cap H)) dH.$$
\endproof
It is explained in \cite{ComteMerle} that $\sigma_k(X,0)=1$ if $0 \le k \le d_0$, so if $k < d_0$ then $\lim_{\epsilon \rightarrow 0} \frac{\Lambda_k(X , X \cap B_\epsilon)}{b_k \epsilon^k}=0$.

\end{document}